\documentclass{amsart}
\usepackage{enumitem}
\usepackage{amscd}
\usepackage{amsfonts}
\usepackage[small,nohug]{diagrams}
\usepackage{bm}
\usepackage{pstricks}
\usepackage{pst-plot}

\newtheorem{theorem}{Theorem}[section]

\newtheorem{proposition}[theorem]{Proposition}

\theoremstyle{remark}

\theoremstyle{definition}
\newtheorem{definition}[theorem]{Definition}
\newtheorem{example}[theorem]{Example}

\numberwithin{equation}{section}

\begin{document}

\title[Illustrating an error in \cite{P}] 
{Illustrating an error in ``An equivalent condition for a uniform space to be coverable''}
\author{Brendon LaBuz}
\address{Saint Francis University, Loretto, PA 15940}
\email{blabuz@@francis.edu}
\date{\today}

\begin{abstract}
Berestovskii and Plaut introduced the concept of a coverable uniform space \cite{BP} when developing their theory of generalized universal covering maps for uniform spaces. Brodskiy, Dydak, LaBuz, and Mitra introduced the concept of a locally uniformly joinable uniform space \cite{Rips} when developing their theory of generalized uniform covering maps which was motivated by \cite{BP}. It is easy to see that a chain connected coverable uniform space is locally uniformly joinable. This paper points out an error in the attempt in \cite{P} to prove that a locally uniformly joinable chain connected uniform space is coverable.
\end{abstract}

\maketitle
\tableofcontents

\medskip
\medskip

\section{Introduction}
Given a uniform space $X$ and an entourage $E$ of $X$, an $E$-chain in $X$ is a sequence $\{x_1,\ldots, x_n\}$ in $X$ such that $(x_i,x_{i+1})\in E$ for each $i<n$. We are interested in defining equivalence classes of $E$-chains. In \cite{BP} this is done by considering $E$-homotopies between $E$-chains. An $E$-homotopy is a sequences of moves where each move consists of adding or deleting a point of the $E$-chain (except either endpoint) so that the result is another $E$-chain. In \cite{Rips} this procedure is framed in terms of Rips complexes; the two concepts are the same and we will use the former in this paper. Following \cite{BP} and \cite{P}, define $X_E$ as follows. Fix a basepoint $*$ of $X$. Then $X_E$ is the set of $E$-homotopy classes $[\alpha]_E$ of $E$-chains $\alpha=\{*=x_0, \ldots,x_n\}$. Give $X_E$ the uniform structure whose basis consists of the sets $F^*$, where for each entourage $F\subset E$ of $X$, $F^*$ consists of ordered pairs $([\alpha]_E,[\beta]_E)$ of elements of $X_E$ where $[\alpha]_E=[x_0,\ldots,x_{n-1},x]_E$ and $[\beta]_E=[x_0,\ldots,x_{n-1},y]_E$ and $(x,y)\in E$. This condition is the same as saying that $\overline\alpha*\beta$ is $E$-homotopic to the $E$-chain $\{x,y\}$ where $x$ and $y$ are the endpoints of $\alpha$ and $\beta$ ($\overline\alpha$ denotes the reverse of $\alpha$).

We will consider the inverse limit of the spaces $X_E$ where the limit is over all entourages of $X$. The bonding maps are the maps $\phi_{EF}:X_F\to X_E$ ($F\subset E$) where an element of $X_F$ is simply considered as an $E$-chain. Let $\phi_{XE}:X_E\to X$ denote the endpoint mapping. The inverse limit of $\{X_E,\phi_{EF}\}$ will be denoted as $\widetilde X$ and is given the inverse limit uniformity where entourages of $\widetilde X$ are $\phi_E^{-1}(F^*)$ where $\phi_E:\widetilde X\to X_E$ is the projection and $F^*$ is an entourage of $X_E$. Here we use the usual convention that for a map $f$ between uniform spaces, $f^{-1}(E)$ means $(f\times f)^{-1}(E)$. Note that in \cite{Rips} elements of the inverse limit are called generalized paths and the limit is denoted as $GP(X,x_0)$. Also, one can consider the set $GP(X)$ of all generalized paths in $X$ without restricting to starting at a particular basepoint. In \cite{Rips} entourages of $GP(X)$ (and also of $GP(X,x_0)$) are denoted as $E^*$ and are defined directly. In the current notation these entourages are $\phi_E^{-1}(E^*)$ and the set of all such entourages forms a basis for the uniform structure on $\widetilde X$ since for $F\subset E$, $\phi_F^{-1}(F^*)\subset \phi_E^{-1}(F^*)$.

A uniform space $X$ is chain connected if for each entourage $E$ of $X$, every pair of points $x,y\in X$ can be joined by an $E$-chain. For chain connected spaces the construction of $X_E$ and therefore $\widetilde X$ is independent of the choice of basepoint \cite[Remark 18]{BP}.

A uniform space $X$ is defined to be locally uniformly joinable \cite{Rips} if for each entourage $E$ of $X$ there is an entourage $F\subset E$ such that such that if $(x,y)\in F$ then $x$ and $y$ can be joined by a generalized path $\alpha$ where the $E$-chain associated with $\alpha$ is $E$-homotopic to the $E$-chain $\{x,y\}$. We call such a generalized path $E$-short. Note that in the initial preprint of \cite{Rips} this concept was referred to as uniformly joinable. The authors subsequently changed the terminology to locally uniformly joinable. Thus references to uniform joinability in \cite{P} really mean local uniform joinability.

In order to give a definition of local uniform joinability that avoids the term ``generalized path," suppose that $X$ is chain connected and define $X$ to be locally uniformly joinable if for every entourage $E$ there is an entourage $F\subset E$ so that given $(x,y)\in F$ there is an element $([\alpha_H])_H\in\widetilde X$ (where the basepoint is taken to be $x$) such that $\alpha_E$ is $E$-homotopic to $\{x,y\}$. Thus we see that if $X$ is locally uniformly joinable chain connected, for every entourage $E$ there exists an entourage $F\subset E$ such that $\phi_{EF}(X_F)\subset\phi_E(\widetilde X)$. That means that the inverse system $\{X_H\}_H$ satisfies the strong Mittag-Leffler condition (\cite[Definition 4.1]{L}). Conversely, assuming $X$ is chain connected, if the inverse system is strong Mittag-Leffler then one can see that $X$ is locally uniformly joinable.

A uniform space $X$ is defined to be coverable \cite{BP} if there is a basis of entourages for $X$ such that for each entourage $E$ in the basis, the projection $\phi_E:\widetilde X\to X_E$ is surjective. Thus we can easily see that for chain connected spaces coverability implies local uniform joinability. Given $E$ one simply takes $F\subset E$ so that $F$ is a member of the basis of ``coverable entourages." One immediately sees that coverability is on the face of it more difficult to satisfy. We are requiring that for every entourage $E$ there exists an entourage $F\subset E$ so that for every $(x,y)\in F$, there is a generalized path joining $x$ and $y$ that is $F$-short.

\section{Uniform openness}

In \cite{P}, a proof is offered that a locally uniformly joinable chain connected space is coverable. Let us see that the hypothesis of Proposition 10 in \cite{P} is equivalent to the space being locally uniformly joinable chain connected. The hypotheses is that $X$ is a chain connected uniform space such that for every entourage $E$ of $X$, the projection $\phi_E:\widetilde X\to X_E$ has uniformly open image in $X_E$. We first give the definition of uniformly open. Recall that given an entourage $E$ of $X$ and a point $x\in X$, $B(x,E)=\{y\in X:(x,y)\in E\}$.

\begin{definition}\cite{P}
A subset $A$ of a uniform space $X$ is uniformly open if there is an entourage $E$ of $X$ such that $B(x,E)\subset A$ for all $x\in A$.
\end{definition}

James defines what it means for a map between uniform spaces to be uniformly open.

\begin{definition} \cite{J}
A function $f:X\to Y$ between uniform spaces is uniformly open if for each entourage $E$ of $X$ there is an entourage $F$ of $Y$ such that $B(f(x),F)\subset f(B(x,E))$ for all $x\in X$.
\end{definition}

There is a nice connection between these two concepts.

\begin{proposition}
A subset $A$ of a uniform space $X$ is uniformly open if and only if the inclusion $A\hookrightarrow X$ is uniformly open.
\end{proposition}

\begin{proof}
Suppose the inclusion map is uniformly open. Then there is an entourage $F$ of $X$ so that $B(x,F)\subset B(x,A\times A)=A$ for all $x\in A$. Now suppose $A$ is uniformly open. Suppose $E$ is an entourage of $X$. Then $E\cap (A\times A)$ is an arbitrary entourage of $A$. Let $F$ be an entourage of $X$ so that $B(x,F)\subset A$ for each $x\in A$. Then $B(x,F\cap E)\subset B(x,E\cap (A\times A))$ for all $x\in A$ so the inclusion is uniformly open.
\end{proof}

\begin{proposition}\label{UnifOpenLocallyUnifJoin}
The following are equivalent for a chain connected uniform space $X$.
\begin{enumerate}[label=\normalfont{\arabic*.}]
\item $X$ is uniformly joinable.
\item For each entourage $E$ of $X$, the image $\phi_E(\widetilde X)$ is uniformly open in $X_E$.
\item The endpoint map $\phi:\widetilde X\to X$ is uniformly open.
\end{enumerate}
\end{proposition}

\begin{proof}
The proposition is proved directly even though some of it could be proven by using results in \cite{Rips} and \cite{P}.

(1.$\implies$2.) Given an entourage $E$ of $X$ choose an entourage $F$ of $X$ so that any $(x,y)\in F$ can be joined by an $E$-short generalized path. Suppose $([\alpha_H]_H)\in\widetilde X$ and $([\alpha_E]_E,[\beta]_E)\in F^*$ for some $[\beta]_E\in X_E$. We show $[\beta]_E\in\phi_E(\widetilde X)$. There is an $E$-short generalized path $([\gamma]_H)_H$ from the endpoint of $\alpha_E$ to the endpoint of $\beta$. Notice $[\alpha_E* \gamma_E]_E=[\beta]_E$.

(2.$\implies$3.) Consider a basic entourage $\phi_E^{-1}(E^*)$ of $\widetilde X$. Choose an entourage $F^*$ of $X_E$ so that for each $([\alpha_H]_H)\in\widetilde X$, $B([\alpha_E]_E,F^*)\subset \phi_E(\widetilde X)$. To see that $\phi:\widetilde X\to X$ is uniformly open, suppose $([\alpha_H]_H)\in\widetilde X$ has endpoint $x$ and suppose $(x,y)\in F$ for some $y\in X$. Let $\beta$ be the $E$-chain obtained by adding $y$ to the end of $\alpha_E$. Then $([\alpha_E]_E,[\beta]_E)\in F^*$ so $[\beta]_E=\phi_E(([\beta_H]_H))$ for some $([\beta_H]_H)\in\widetilde X$. Note $(([\alpha_H]_H),([\beta_H]_H))\in\phi_E^{-1}(F^*)\subset \phi_E^{-1}(E^*)$.

(3.$\implies$1.) Given an entourage $E$ of $X$, choose an entourage $F$ of $X$ so that $B(\phi(([\alpha_H]_H)),F)\subset \phi(B(([\alpha_H]_H), \phi_E^{-1}(E^*)))$ for any $([\alpha_H]_H) \in \widetilde X$. We show $\phi_{EF}(X_F)\subset \phi_E(\widetilde X)$. We will use induction on lengths of chains. The base case is taken care of by the constant generalized path. Suppose an $F$-chain $\alpha=\{x_0,\ldots,x_n\}$ has $[\alpha]_E=\phi_E(([\alpha_H]_H))$ for some $([\alpha_H]_H)\in\widetilde X$ and $(x_n,x_{n+1})\in F$. Then by hypothesis $x_{n+1}$ is the endpoint of some $([\beta_H]_H)\in \widetilde X$ such that $(([\alpha_H]_H),([\beta_H]_H))\in \phi_E^{-1}(E^*)$. But then $[\beta_E]_E=[x_0,\ldots,x_{n+1}]_E$.
\end{proof}

\section{The error in \cite{P}}

This section contains an example that shows the proofs of both Propositions 9 and 10 in \cite{P} are not valid and points out an error in the proof of Proposition 9.

Consider Proposition 10 in \cite{P}. The hypothesis is that $X$ is chain connected and $\phi_E(\widetilde X)$ is uniformly open in $X_E$ for each entourage $E$ of $X$ so by Proposition \ref{UnifOpenLocallyUnifJoin} any locally uniformly joinable chain connected space (in particular a path connected and uniformly locally path connected space) will satisfy the hypothesis. 

Now consider the proof of Proposition 10 in \cite{P}. First note that for a uniform space $X$ and an entourage $E$ of $X$, $\phi_{XE}:X_E\to X$ is the notation for the endpoint mapping. In this proof, we are considering the entourage $G=\phi_E^{-1}(E^*)$ of $\widetilde X$ and the endpoint mapping $\phi_{\widetilde X G}:\widetilde X_G\to \widetilde X$ where the basepoint of $\widetilde X_G$ is taken to be the constant generalized path. Indeed, the beginning of the proof of Proposition 10 considers an arbitrary entourage $E$ of $X$ and refers to the diagram (1) in the statement of Proposition 9 where $G$ is defined as $\phi_E^{-1}(E^*)$. 

Toward the end of the proof of Proposition 10 it is claimed that $\phi_{\widetilde X G}$ is a uniform homeomorphism which means in particular that it is a bijection. The following is an example of a path connected and uniformly locally path connected space $X$ and an entourage $E$ of $X$ for which $\phi_{\widetilde X G}$ is not injective, showing that the proof of Proposition 10 is not valid. 

\begin{figure}\label{hexagon}
\begin{pspicture*}(-1,-1)(7,5)
\psline(0,1)(2,0)
\psline(2,0)(4,0)
\psline(4,0)(6,1)
\psline(6,1)(4,2)
\psline(4,2)(2,2)
\psline(2,2)(0,1)
\psline(0,1)(0,4)
\psline(0,4)(3,4)
\psline(3,4)(3,1)
\rput(-0.3,1){$a$}
\rput(2,-0.3){$b$}
\rput(4,-0.3){$c$}
\rput(6.3,1){$d$}
\rput(4,2.3){$e$}
\rput(2,2.3){$f$}
\rput(3,0.7){$o$}
\rput(-0.3,4){$g$}
\rput(3.3,4){$h$}
\end{pspicture*}
\caption{The map $\phi_{\widetilde X G}$ is not injective for this space and a particular entourage $E$.}
\end{figure}

\begin{example} \label{counterexample}
Consider a regular hexagon whose sides have length 1. Label its vertices in order from $a$ to $f$. Add the center $o$ of the hexagon and a vertical square with base $ao$ that we remove. Label the vertices of the square in order as $a, g, h, o$ (see Figure \ref{hexagon}). Consider the entourage $E=\{(x,y):d(x,y)\leq 1\}$ of the resulting space $X$ where the metric is inherited from the standard metric on $\mathbb R^3$. Let $a=*$ be the basepoint of $X$. Then the basepoint of $\widetilde X$ is the constant element $\iota=([*]_H)$. Consider the $G=\phi_E^{-1}(E^*)$-chain $\alpha=\{\iota,ab,abc,abcd,afed,afe,af,\iota\}$ in $\widetilde X$ where each element of the chain is taken to be the element of $\widetilde X$ induced by the indicated path in $X$. Now $\phi_{\widetilde X G}([\alpha]_G)=\phi_{\widetilde X G}([\iota]_G)=\iota$. We can see that $[\alpha]_G\neq[\iota]_G$. Indeed, any $\phi_E^{-1}(E^*)$-homotopy from $\alpha$ to $\iota$ would need to introduce an element of $\widetilde X$ whose endpoint is $o$. But the only such element is the one $\beta$ that is induced by the path $agho$ and $\beta$ is not $\phi_E^{-1}(E^*)$-close to any element of $\widetilde X$ with endpoint on the hexagon since $\phi_E(\beta)=[a,g,h,o]_E$ which is not the same equivalence class as $[a,o]_E$.
\end{example}

Notice the similarity to the examples contained in Section 7 of \cite{Rips} which show that certain strategies for showing a locally uniformly joinable space is coverable cannot work.

Next consider Proposition 9 in \cite{P}. The hypothesis is the following. Suppose $X$ is chain connected and $E$ is an entourage in $X$ such that $A=\phi_E(\widetilde X)$ is uniformly open in $X_E$. Set $A=\phi_E(\widetilde X)$ and $D=E^*\cap (A\times A)$ which is an entourage of $A$. Then we set $G=\phi_E^{-1}(D)=\phi_E^{-1}(E^*)$. Now $A_D$ is the space of $D$-chains in $A$ and there is a map $\theta:\widetilde X_G\to A_D$ induced by $\phi_E$. Then the conclusion (in part) is that there is a map $\psi:\widetilde X\to A_D$ that makes the following diagram commute.

\begin{diagram}
\widetilde X_G & \rTo^{\phi_{\widetilde X G}}  & \widetilde X \\
\dTo^{\theta}  & \ldTo<\psi      &           \\
A_D							&									&						\\
\end{diagram}

We will describe the construction of $\psi$ and explain the error in the argument that the diagram commutes.

First, an entourage $E_A=\phi_{XE}(D)$ of $X$ is defined. Since $A$ is uniformly open, $E_A$ is an entourage of $X$ (\cite[Lemma 5]{P}). Then, given $([\alpha_E]_E) \in \widetilde X$, $\psi(([\alpha_E]_E))$ is defined to be the equivalence class  $[[x_0]_E,[x_0,x_1]_E,\ldots ,[x_0,x_1,\ldots ,x_n]_E]_D$, where $[x_0,\ldots,x_n]_{E_A}=[\alpha_{E_A}]_{E_A}$.

Following the proof of Proposition 9, suppose $\eta=\{*=y_0,y_1,\ldots ,y_n\}$ is an $\phi_E^{-1}(E^*)$-chain in $\widetilde X$. For $i\leq n$, let $x_i$ be the endpoint of $y_i$. Then  $\{x_0,\ldots,x_n\}$ is an an $E_A$-chain. Then it is claimed that $[[*]_E,[*,x_1]_E,\ldots ,[*,x_1,\ldots ,x_n]_E]_D=\psi\circ\phi_{\widetilde X G}([\eta]_G)$ which is incorrect. Even though $\{x_0,\ldots,x_n\}$ is an $E_A$-chain, we cannot know that $[x_0,\ldots,x_n]_{E_A}=(y_n)_{E_A}$. The same Example \ref{counterexample} illustrates this fact.

Let us follow the proof of Proposition 9 using the same $X$ and $E$ from Example \ref{counterexample}. Let $\eta =\{*=\iota,ab,abc,abcd,afed,afe,af,\iota\}$, the same $\phi_E^{-1}(E^*)$-chain used in Example \ref{counterexample}. Now $\psi \circ \phi_{\widetilde X G}([\eta]_{G})$ is the equivalence class of the constant chain, but $\theta([\eta]_{G})=[[a]_E,[a,b]_E,[a,b,c]_E,[a,b,c,d]_E,[a,f,e,d]_E,[a,f,e]_E,[a,f]_E,[a]_E]_D$ is not the equivalence class of the constant chain. There is an $E^*$-homotopy in $X_E$ from $\{[a]_E,[a,b]_E,[a,b,c]_E,[a,b,c,d]_E,[a,f,e,d]_E,[a,f,e]_E,[a,f]_E,[a]_E\}$ to the constant chain, but not a $D$-homotopy in $A$. Thus the diagram does not commute. In fact, this example is a counterexample for Proposition 9 for this particular $\psi$.

\end{document}